\documentclass[12pt]{article}
\usepackage{graphicx,latexsym,amssymb,amsmath,amsthm}

\begin{document}

\baselineskip=18pt
\newtheorem{tetel}{Theorem}[section]
\newtheorem{dfn}[tetel]{Definition}
\newtheorem{all}[tetel]{Claim}
\newtheorem{lem}[tetel]{Lemma}
\newtheorem{kov}[tetel]{Corollary}
\newcommand{\bp}{\noindent\emph{Proof. }}
\newcommand{\ep}{\hfill$\Box$\par\medskip}
\newcommand{\nop}{\hfill$\Box$\par}
\newcommand{\megj}{\noindent{\bf Remark. }}

\title{The ultimate question}
\author{G\'abor Wiener \\ Department of Computer Science and Information Theory \\ Budapest University of Technology and Economics \\ \\ Makoto Araya \\ Department of Computer Science \\ Shizuoka University} 
\maketitle

\begin{abstract}
We present a planar hypohamiltonian graph on 42 vertices and show some consequences.
\end{abstract}

\section{Introduction and results}

A graph is called hypohamiltonian if it is not Hamiltonian but deleting any vertex gives a Hamiltonian graph. Hypohamiltonian graphs were extensively studied by Sousselier, Thomassen, Chvatal, Hatzel, Zamfirescu, Aldred et al., and many others, see for example the survey by Holton and Sheehan \cite{HS}. Chvatal raised the question in 1973 whether there exists a planar hypohamiltonian graph. This was answered by Thomassen, who found such a graph on 105 vertices in 1976. Hatzel (in 1979) found a smaller planar hypohamiltonian graph, having 57 vertices and this was improved to 48 vertices by Zamfirescu and Zamfirescu in 2007.  Here we present a planar hypohamiltonian graph on 42 vertices, making possible that the Ultimate Question in Douglas Adams's The Hitch Hikers' Guide to the Galaxy was "What is the size of the smallest planar hypohamiltonian graph?". 

Using a theorem of Thomassen we can also construct a planar hypotraceable graph (a graph is hypotraceable if it does not contain a Hamiltonian path but the deletion of any vertex gives a graph having a Hamiltonian path) on 162 vertices improving the best known bound of 186. 

Denote the smallest number of vertices of a planar $k$-connected graph, in which every $j$ vertices are omitted by some longest cycle (path) by $\overline{C_k^j}$ ($\overline{P_k^j}$). The bounds on the numbers $\overline{C_3^1}$, $\overline{C_3^2}$, $\overline{P_3^1}$, $\overline{P_3^1}$ are also improved using our construction. 

\section{The construction}

Our graph is the following.

\bigskip

\includegraphics[height=5cm]{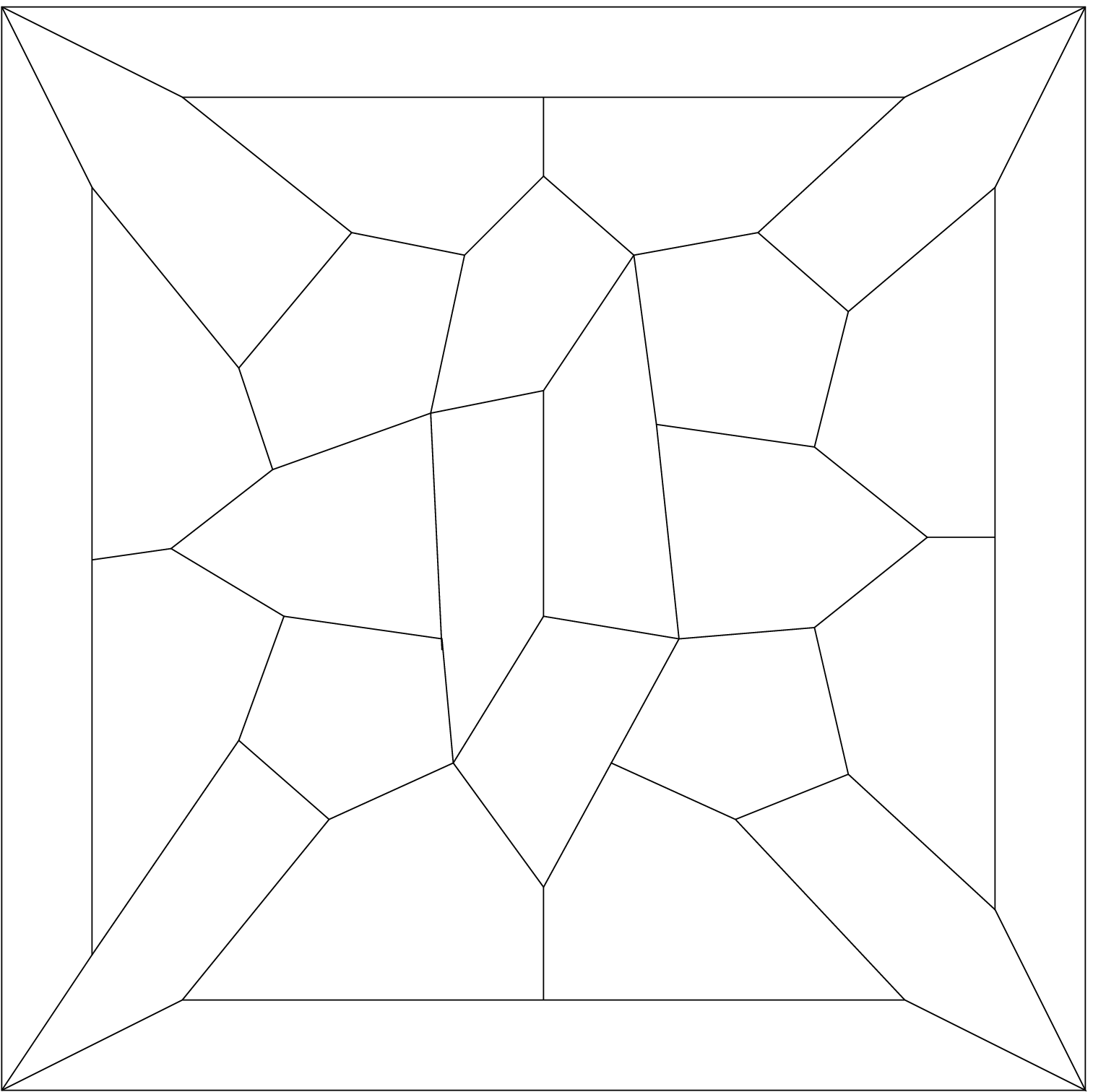}

\medskip

It is not Hamiltonian by Grinberg's criterion, and using i.e. Mathematica it can be easily verified that the deletion of any vertex gives a Hamiltonian graph. 
%\section{}

%\section{Introduction}

%\begin{lem}

%\end{lem}

%\begin{lem}
 
%\end{lem}

\end{document}